\documentclass[11pt,a4paper,twoside,leqno]{article}
\usepackage{amssymb}

\newcounter{rema}

\newcommand{\RR}{\mathbb{R}}

\newcommand{\rf}[1]{(\ref{#1})}
\newtheorem{theorem}{Theorem}

\newtheorem{theof}{Th\'eor\`eme}

\newcommand{\sectio}[1]{
              \par\bigskip
              \stepcounter{section}
\noindent
{\bf \arabic{section}. #1}\par \bigskip}

\begin{document}

\noindent
{\bf \today}

\noindent
\'Equations aux d\'eriv\'ees partielles/Partial Differential Equations

\vskip1.5cm

\noindent
{\Large\bf  Asymptotic profiles of solutions to convection-diffusion
equations}
\bigskip

\noindent 
{\bf Sa\"\i d BENACHOUR,$^1$ Grzegorz KARCH,$^2$ and Philippe
  LAUREN\c COT$^3$} \\ 
{
$^1$ Institut Elie Cartan-Nancy, Universit\'e Henri Poincar\'e, 
 BP 239, F-54506 Vand\oe uvre les Nancy cedex, France. 
E-mail: {\it benachou@iecn.u-nancy.fr} \\
$^2$ Instytut Matematyczny, Uniwersytet Wroc\l awski, 
pl. Grunwaldzki 2/4, 50-384 Wroc\l aw, Poland, and 
Institute of Mathematics, Polish Academy  of Sciences, Warsaw
(2002-2003). E-mail: {\it karch@math.uni.wroc.pl}\\
$^3$ Math\'ematiques pour l'Industrie et la Physique, 
CNRS UMR 5640, Universit\'e Paul Sabatier-Toulouse 3,
118 route de Narbonne, 31062 Toulouse cedex 4, France. E-mail: {\it
  laurenco@mip.ups-tlse.fr} 
}

\bigskip

\hrule

\medskip

{\bf Abstract.}
The large time behavior of zero-mass solutions to the Cauchy problem
for the convection-diffusion equation $u_t- u_{xx}+(|u|^q)_x=0,$
$u(x,0)=u_0(x)$ is studied when $q>1$ and the initial datum $u_0$
belongs to $L^1(\RR, (1+|x|)\;dx)$ and satisfies $\int_\RR
u_0(x)\ dx=0$. We provide conditions on the size, and shape of the
initial datum $u_0$ as well as on the exponent $q>1$ such that the
large time asymptotics of solutions is given either by the derivative
of the Gauss-Weierstrass kernel, or by a self-similar solution of the
equation, or by hyperbolic $N$-waves.

\medskip

\hrule
 
\medskip

\noindent
{\large\bf  Comportement asymptotique des solutions d'\'equations de
convection-diffusion} 

\medskip

{\bf R\'esum\'e.} Le comportement asymptotique des solutions de masse
nulle du probl\`eme de Cauchy pour l'\'equation de
convection-diffusion $u_t- u_{xx}+(|u|^q)_x=0,$ $u(x,0)=u_0(x)$ est
\'etudi\'e lorsque $q>1$ et la donn\'ee initiale $u_0$ appartient \`a
$L^1(\RR, (1+|x|)\;dx)$ et satisfait $\int_\RR u_0(x)\ dx=0$. Nous
donnons des conditions sur l'amplitude et la forme de la donn\'ee
initiale $u_0$ et sur l'exposant $q>1$ sous lesquelles le comportement
asymptotique des solutions est d\'ecrit par la d\'eriv\'ee premi\`ere
du noyau de Gauss-Weierstrass, ou par une solution auto-similaire de
l'\'equation, ou par une $N$-onde hyperbolique. 

\medskip

\hrule
 
\medskip

\noindent
{\bf Version fran\c caise abr\'eg\'ee} 
\bigskip

Nous \'etudions le comportement asymptotique, lorsque $t\to\infty$,
des solutions de masse nulle du probl\`eme de Cauchy pour l'\'equation
de convection-diffusion \rf{eq}-\rf{ini}. Il est clair que cette
\'equation met en jeu une comp\'etition entre le terme de diffusion
$u_{xx}$ et le terme de convection non lin\'eaire $(\vert
u\vert^q)_x$, et le but de notre \'etude est de d\'eterminer lequel de
ces deux termes devient pr\'epond\'erant lorsque $t\to\infty$. Lorsque la
masse totale $M=\int_\RR u_0(x)\ dx \ne 0$, ce probl\`eme a \'et\'e
consid\'er\'e en particulier dans \cite{EZ,EVZ1,kim}. Lorsque $M=0$, les
r\'esultats de \cite{EZ,EVZ1,kim} se traduisent seulement par la
convergence de $u(t)$ vers z\'ero dans $L^1(\RR)$. Nous pr\'esentons,
dans cette note, des r\'esultats plus pr\'ecis pour certaines classes
de donn\'ees initiales. En effet, nous exhibons le premier terme dans le
d\'eveloppement asymptotique des solutions de \rf{eq}-\rf{ini} selon
les valeurs du param\`etre $q>1$ et le signe (et \'eventuellement la
taille) de la primitive de la donn\'ee initiale $u_0$.

Nous supposons que $u_0$ v\'erifie \rf{as1} et
d\'esignons par $u$ la solution correspondante de
\rf{eq}-\rf{ini}. Signalons d\`es \`a pr\'esent que notre approche est
bas\'ee sur l'\'etude de l'\'equation de Hamilton-Jacobi diffusive
\rf{eq+}-\rf{ini+} qui est naturellement associ\'ee \`a
\rf{eq}-\rf{ini}. En effet, la primitive $U$ de $u$, d\'efinie par \rf{U},
est solution de \rf{eq+}-\rf{ini+}.
 
En premier lieu, nous \'etudions le cas o\`u la diffusion r\'egit le
comportement en temps grands.

\begin{theof} \label{th:a}
Supposons que l'une des trois conditions suivantes soit v\'erifi\'ee~: 
(i) $U_0\ge 0$ et  $q>{3/2}$, (ii) $U_0\le 0$ et  $q\ge 2$, ou (iii)
$U_0\le 0$,  $q\in (3/2, 2)$ et $u_0$ satisfait \rf{small:u0}. Alors,
$u\in L^q(\RR\times (0,\infty))$ et 
$$
\lim_{t\to\infty} t^{(1-1/p)/2+1/2}\  \|u(t)-I_\infty\ G_x(t)\|_p=0\,,
\quad p\in [1,\infty]\,, 
$$
o\`u $I_\infty$ est un r\'eel non nul d\'efini par \rf{I:infty} et $G$
est la solution fondamentale de l'\'equation de la chaleur.
\end{theof}

Consid\'erons ensuite le cas o\`u $U_0\ge 0$ et $q\in
(1,3/2)$. L'influence des effets convectifs et diffusifs s'\'equilibre
lorsque $t\to\infty$ et on a le r\'esultat suivant~:

\begin{theof}\label{th:b}
Si $q\in (1, 3/2)$ et $U_0\ge 0$ v\'erifie la condition de
croissance \rf{a-u0}, alors 
$$
\lim_{t\to\infty} t^{(1-1/p)/2+a/2} \| u(t)- W_x(t)\|_p=0\,, \quad
p\in [1,\infty]\,, 
$$
o\`u $a:=(2-q)/(q-1)$ et $W$ d\'esigne la \textit{solution tr\`es
  singuli\`ere} de \rf{eq+}. 
\end{theof}

Enfin, nous identifions une classe de donn\'ees initiales $u_0$ pour
lesquelles le comportement en temps grands est r\'egi par la
convection. 

\begin{theof} \label{th:c}
Supposons que $q\in (1,2)$ et que $U_0\le 0$ v\'erifie
\rf{large:u0}. Alors
$$
\lim_{t\to\infty} t^{(1-1/p)/q}\ \|u(t)-N_{\sigma,\sigma}(t)\|_p=0\,,
\quad p\in [1,\infty)\,,
$$
o\`u $\sigma$ est le r\'eel strictement positif d\'efini par
\rf{evian} et $N_{\sigma,\sigma}$ est la $N$-onde hyperbolique
d\'efinie par \rf{N-wave}. 
\end{theof}

Lorsque $q\in(1,4/(1+\sqrt{3}))$, le Th\'eor\`eme~\ref{th:c} est vrai sans la
condition de taille \rf{large:u0} sur $u_0$ et nous conjecturons qu'il
reste valable sans cette condition pour tout $q\in (1,3/2)$. 

On notera que, lorsque $q\in (3/2,2)$ et selon la taille de la
donn\'ee initiale, on peut observer aussi bien un comportement
asymptotique o\`u la diffusion domine qu'un comportement o\`u la
convection est pr\'epond\'erante. De plus, il existe des solutions
pour lesquelles les deux effets se compensent \cite{BSW02}.
  
Les preuves des r\'esultats pr\'esent\'es ici, ainsi que des
extensions en dimension sup\'erieure, font l'objet de l'article
\cite{BKL}. 

\medskip

\hrule

\medskip

\sectio{Introduction}

The large time behavior of zero-mass solutions to the Cauchy problem 
\begin{eqnarray}
&u_t- u_{xx}+(|u|^q)_x=0,& x\in \RR\,, \;t>0\,,
\label{eq}\\
&u(x,0)=u_0(x),& x\in\RR\,,
\label{ini}
\end{eqnarray}
is investigated when $q>1$, under the main assumption
\begin{equation}
u_0\in L^1(\RR, (1+|x|)\;dx)\cap W^{1,\infty}(\RR)\,, \quad
  u_0\not\equiv 0 \quad \mbox{and}\quad \int_\RR u_0(x)\ dx=0\,. 
\label{as1}
\end{equation} 
The equation (\ref{eq}) includes two competing effects, namely the
diffusion $u_{xx}$ and the nonlinear convection $(|u|^q)_x$, and the
main issue in the study of the large time asymptotics
is to figure out whether one of these effects dominates for large
times. For integrable initial data with non-zero mass, i.e., 
$$
u_0\in L^1(\RR) \quad \mbox{and}\quad \int_\RR u_0(x)\ dx=M\ne 0\,,
$$
this question has already been investigated and the outcome may be
summarized as follows: if $q>2$, the large time dynamics is dominated
by the diffusion and $u$ behaves as the fundamental (or source-type)
solution $M\, G$ to the linear heat equation with $G(x,t)=(4\pi
t)^{-n/2}\exp(-|x|^2/(4t))$ for $(x,t)\in\RR\times (0,\infty)$, see,
e.g., \cite{EZ}. If $q\in (1,2)$ and either $u_0\ge 0$ (or $u_0\le 0$)
\cite{EVZ1} or 
\begin{equation}
\label{kimas}
\inf_{x\in\RR} \int_{-\infty}^x u_0(y)\ dy \le 0
\end{equation} 
\cite{kim}, the convection term is preponderant for large times and
there is $\sigma\in\RR$ such that $u$ behaves as the $N$-wave solution
$N_{\sigma,\sigma+M}$ to the nonlinear conservation law $z_t+(\vert
z\vert^q)_x=0$. Moreover, $\sigma=0$ if $u_0\ge 0$, $\sigma=-M$ if
$u_0\le 0$ and $\sigma\ge 0$ under the assumption
(\ref{kimas}). Recall that the $N$-wave solution $N_{\alpha,\beta}$ is
given explicitly by  
\begin{equation}
N_{\alpha,\beta}(x,t)= 
\left\{
\begin{array}{ccc}
\displaystyle{ {\rm sign}\, x\cdot \left( {|x|\over qt}
  \right)^{1/(q-1)} } , & & 
\displaystyle{-q\left({\alpha\over q-1}\right)^{(q-1)/q} \leq
  {x\over t^{1/q}}
\leq q\left({\beta\over q-1}\right)^{1/(q-1)}\,,}
\\
& & \\
0 & &\mbox{otherwise\,,}
\end{array}
\right.
\label{N-wave}
\end{equation}
for $\alpha\ge 0$ and $\beta\ge 0$ (see, e.g., \cite{Sm}). When $q=2$,
there is a balance between the diffusive and convective terms and $u$
behaves as the unique self-similar source-type solution $S_M$ to
(\ref{eq}), that is, $S_M$ is the unique solution to (\ref{eq}) with
initial datum $M\, \delta$ \cite{EZ}. 

The previous results remain of course valid when $M=0$. However, when
$q\ge 2$, they reduce to the convergence of $u(t)$ to zero in
$L^1(\RR)$, while if $q\in (1,2)$ and (\ref{kimas}) is fulfilled, the
possibility that $\sigma=0$ is not excluded in \cite{kim}. 
Obtaining non-vanishing intermediate asymptotics in the zero-mass case
thus seems to be more delicate, and is the purpose of this note. On
the one hand, we identify conditions under which the constant $\sigma$
defined previously does not vanish, thus leading to a truly
hyperbolic large time behavior. On the other hand, we complete the
study performed in \cite{KS} where conditions on the initial data
$u_0$ are found for which the diffusion still plays a role in the
large time asymptotics. Indeed, under the assumption that there is
$\beta\in (0,1)$ 
such that $\sup_{t\ge 0} t^{\beta/2}\ \|e^{t\Delta}u_0\|_1<\infty$ if
$q\ge 2$ and $\sup_{t\ge 0} t^{\beta/2}\ \|e^{t\Delta}u_0\|_1$ is
sufficiently small for $q\in (1+1/(1+\beta),2)$, the large time
asymptotics is dominated by the diffusive term and the asymptotics of
solutions to \rf{eq}-\rf{ini} is described by self-similar solutions
to the heat equation \cite{KS}. A balance between diffusion and
convection is also observed in \cite{KS} in the critical case
$q=1+1/(1+\beta)$ and, for suitably small initial data, the
asymptotics of solutions to \rf{eq}-\rf{ini} is described by a new
class of self-similar solutions to \rf{eq}. 

\medskip

>From now on, we assume that $u_0$ fulfils (\ref{as1}). We denote by
$u=u(x,t)$ the corresponding solution to \rf{eq}-\rf{ini} and set
\begin{equation}
U(x,t) =\int_{-\infty}^x u(y,t)\;dy = -\int_x^\infty u(y,t)\;dy\,,
\label{U}
\end{equation}
the second inequality being true since the integral of $u(t)$ over
$\RR$ vanishes for each $t\ge 0$ by (\ref{as1}). Owing to (\ref{as1}),
it is easy to see that $U_0:=U(0)$ belongs to $L^1(\RR)\cap
L^\infty(\RR)$ with
$$
\int_\RR U_0(x)\;dx =-\int_\RR xu_0(x)\;dx\,,
$$
and it readily follows from \rf{eq} that $U$ is a solution to the
viscous Hamilton-Jacobi equation 
\begin{eqnarray}
&&U_t-U_{xx} +|U_x|^q=0, \quad x\in\RR,\; t>0, \label{eq+}\\
&&U(x,0)=U_0(x),\quad x\in \RR\,. \label{ini+}
\end{eqnarray}
This observation is the basis of our approach towards the study of the
large time behavior of zero-mass solutions to \rf{eq}-\rf{ini}, since
we will actually investigate the large time dynamics of $U$ and then deduce
corresponding results for $u$.

%%%%%%%%%%%%%%%%%%%%%%%%%%%%%%%%%%%%%%%%%%%%%%%%%%%%%%%%%%%%%%%%%%%%%%%

\sectio{ Diffusion-dominated case}

We begin with the case when the diffusive term is preponderant for
large times.

\begin{theorem} \label{th:01}
Assume that $u_0$ fulfils \rf{as1} and that one of the following
assertions hold true:

(i) $U_0\ge 0$ and  $q>{3/2}$,

(ii) $U_0\le 0$ and  $q\geq 2$,

(iii) $U_0\le 0$,  $q\in (3/2, 2)$, and
\begin{equation}
\left|\int_{\RR}xu_0(x)\;dx\right| \|u_0\|_\infty^{2q-3}
\quad \mbox{is sufficiently small.} \label{small:u0}
\end{equation}

Then, $u\in L^q(\RR\times (0,\infty))$, the constant $I_\infty$ defined
by 
\begin{equation}
I_\infty := - \lim_{t\to\infty} \int_{\RR}xu(x,t)\;dx=
-\int_{\RR}xu_0(x)\;dx-\int_0^\infty \int_{\RR} |
u(x,s)|^q\;dx\,ds
\label{I:infty}
\end{equation}
is well-defined with $I_\infty>0$ if $U_0\ge 0$ and $I_\infty<0$ if
$U_0\le 0$, and 
\begin{equation}
\lim_{t\to\infty} t^{(1-1/p)/2+1/2}\  \|u(t)-I_\infty\
G_x(t)\|_p=0 \;\;\mbox{ for every }\;\; p\in
[1,\infty]\,. \label{lim:G1}  
\end{equation}
\end{theorem}

A result similar to Theorem~\ref{th:01} is also valid for $U$
\cite[Theorems~2.1 \&~2.3]{BKL}. The proof of Theorem~\ref{th:01}
splits in two steps: first, the fact that $u\in L^q(\RR\times (0,\infty))$ and
$I_\infty\in\RR\setminus\{0\}$ under the assumptions of
Theorem~\ref{th:01} follows from \cite{BL99,BK99} in the case (i) and
\cite{LS03} in the case (ii) and (iii). Next, the proof of \rf{lim:G1}
and its analogue for $U$ relies on the representation of $U$ by the
Duhamel formula and $L^\infty$- and $L^q$-estimates on $\partial_x U=u$
established in \cite{BL99,GGK03} and \cite{LS03}, respectively, from
which the estimates
$$
\sup_{t>0}\ \left\{ t^{1/2}\ \Vert u(t)\Vert_1 \right\} + \sup_{t>0}\
\left\{ t\ \Vert u(t)\Vert_\infty \right\} < \infty
$$
follow. Classical properties of the heat semi-group $e^{t\Delta}$ then
allow us to complete the proof of Theorem~\ref{th:01}.
 
%%%%%%%%%%%%%%%%%%%%%%%%%%%%%%%%%%%%%%%%%%%%%%%%%%%%%%%%%%%%%%%%%%%%%%%
%%%%%%%%%%%%%%%%%%%%%%%%%%%%%%%%%%%%%%%%%%%%%%%%%%%%%%%%%%%%%%%%%%%%%%%%

\sectio{Convergence towards very singular solutions}

Our next theorem is devoted to the case $q\in (1, 3/2)$ and $U_0\ge 0$
where there is a balance between the diffusive and convective effects,
and a particular self-similar solution of \rf{eq} appears in
the large time asymptotics.

\begin{theorem}\label{th:02}
Suppose that $u_0$ fulfils \rf{as1}. If $q\in (1, 3/2)$ and $U_0$ is
nonnegative and satisfies
\begin{equation}
\lim_{|x|\to\infty} |x|^a\ U_0(x)=0 \quad \mbox{with}\quad a={2-q\over
q-1}\,, \label{a-u0}
\end{equation}
then
\begin{equation}
\lim_{t\to\infty} t^{(1-1/p)/2+a/2} \| u(t)- W_x(t)\|_p=0
\label{lim:W2}
\end{equation}
for every $p\in [1,\infty]$, where $W$ is the very singular
solution to \rf{eq+} and enjoys the self-similarity property
$W(x,t)=t^{-a/2} W(x t^{-1/2},1)$ for $(x,t)\in\RR\times (0,\infty)$. 
\end{theorem}

Note that $a>1$ for $q<3/2$, hence the convergence rate towards the
self-similar profile in \rf{lim:W2} is faster than the rate in
\rf{lim:G1}. 

We recall that a very singular solution $W$ to \rf{eq+} is a classical
solution to \rf{eq+} in $\RR\times (0,\infty)$ which has a singular
behavior as $t\to 0$, namely
$$
\lim_{t\to 0} \int_{\{\vert x\vert \ge r\}} W(x,t)\ dx = 0 \;\;\mbox{
  and }\;\; \lim_{t\to 0} \int_{\{\vert x\vert \le r\}} W(x,t)\ dx = \infty
$$
for each $r>0$. The existence and uniqueness of the very singular
solution to \rf{eq+} are established in \cite{BKLxx,BL01}. The proof of
Theorem~\ref{th:02} relies on a rescaling method together with some
estimates from \cite{BL01}. More precisely, introducing the sequence
of rescaled functions $U_k(x,t) := k^a\ U(kx,k^2t)$, $k\ge 1$, an easy
computation shows that $U_k$ is a solution to \rf{eq+} while it
follows from \cite{BL01} that
$$
\sup_{t>0}\ \left\{ t^{(a-1)/2}\ \Vert u(t)\Vert_1 + t^{a/2}\ \Vert
  U(t)\Vert_\infty + t^{(a+1)/2}\ \Vert u(t)\Vert_\infty \right\} <
\infty\,. 
$$
Owing to these estimates, one can show that a subsequence of $(U_k)$
converges towards $W$ in $L^p(\RR)$ for every $p\in [1,\infty]$. The
uniqueness of the very singular solution to \rf{eq+} and the use of
the Duhamel representation formula allow to extend the previous
convergence to the whole sequence $(U_k)$ and also to its first
derivative, thus completing the proof of Theorem~\ref{th:02}. 

%%%%%%%%%%%%%%%%%%%%%%%%%%%%%%%%%%%%%%%%%%%%%%%%%%%%%%%%%%%%%%%%%%%%%%%
%%%%%%%%%%%%%%%%%%%%%%%%%%%%%%%%%%%%%%%%%%%%%%%%%%%%%%%%%%%%%%%%%%%%%%%%

\sectio{Hyperbolic-dominated case}

We end up with the case when the large time behavior is ruled by
the convection term and is given by self-similar solutions to the
nonlinear conservation law $z_t+(|z|^q)_x=0$.

\begin{theorem} \label{th:hyp}
Consider $q\in (1,2)$. Assume that $u_0$ fulfils \rf{as1} and that
$U_0$ is non-positive and such that 
\begin{equation}
\|U_0\|_\infty \ \Vert u_{0,x}\Vert_\infty^{1-2/q} \quad \mbox{is
  sufficiently large.} \label{large:u0} 
\end{equation}
Then 
\begin{equation}
\sigma:= \lim_{t\to\infty} \Vert u(t)\Vert_\infty \in (0,\infty)
\label{evian}
\end{equation}
and
\begin{equation}
\lim_{t\to\infty} t^{(1-1/p)/q}\ \|u(t)-N_{\sigma,\sigma}(t)\|_p=0
\label{volvic}
\end{equation}
for every $p\in [1,\infty)$, $N_{\sigma,\sigma}$ being the $N$-wave
solution to $z_t+(|z|^q)_x=0$ defined by \rf{N-wave}. 

In addition, if $q\in (1,4/(1+\sqrt{3}))$, the previous result is valid
without the condition \rf{large:u0} on $u_0$.
\end{theorem}

As a consequence of Theorems~\ref{th:01} and~\ref{th:hyp}, we realize
that, when $q\in (3/2,2)$, both diffusion-dominated and
hyperbolic-dominated large time behaviors are possible, according to
the size of $u_0$. In addition, when $q\in (3/2,2)$, there are
non-positive initial data for which the large time behavior is
governed by a balance beween the diffusive and convective effects:
indeed, self-similar solutions to (\ref{eq+}) of the form
$t^{-(2-q)/(2(q-1))} V(\vert x\vert t^{-1/2})$ are constructed in
\cite{BSW02} and $V$ and its derivatives are smooth and decay
exponentially as $\vert x\vert\to\infty$. 

The main step of the proof of Theorem~\ref{th:hyp} is to show
(\ref{evian}). For that purpose, we use the unilateral estimate
\begin{equation}
U_{xx}(x,t) = u_x(x,t) \le \min{\left\{ \Vert u_{0,x}\Vert_\infty
    \,,\, C(q)\ \Vert U_0\Vert_\infty^{(2-q)/q}\ t^{-2/q}
    \right\}}\,. \label{badoit}  
\end{equation}
This inequality can be seen as a weak form of the Oleinik gradient
estimate for scalar conservation laws and is also derived in
\cite{kim} (see also \cite{FL} for related estimates). 
We however establish a version of (\ref{badoit}) in
a multi-dimensional setting involving the Hessian matrix of $U$
\cite{BKL}. The proof of \rf{badoit} follows from \rf{eq} by the
maximum principle and the $L^\infty$-estimate $\Vert
u(t)\Vert_\infty\le C(q)\ \Vert U_0\Vert_\infty^{1/q}\ t^{-1/q}$
\cite{BL99,GGK03}. Once \rf{evian} is established, the proof of \rf{volvic}
is performed by a classical rescaling method \cite{BKL,kim}.

\begin{remark} (a) Though the conditions \rf{small:u0} and \rf{large:u0} do not
involve the same quantities, it can be shown by an application of
Gagliardo-Nirenberg inequalities that, if \rf{small:u0} is fulfilled,
the quantity involved in \rf{large:u0} is also small.\\
(b) We believe that the condition \rf{large:u0} is not necessary when
$q\in (1,3/2)$ for Theorem~\ref{th:hyp} to hold true but have only
been able to prove it for $1<q<4/(1+\sqrt{3})<3/2$.
\end{remark}

\medskip

The proofs of the results presented in this note, together with their
extension to higher space dimensions are to be found in \cite{BKL}.

%%%%%%%%%%%%%%%%%%%%%%%%%%%%%%%%%%%%%%%%%%%%%%%%%%%%%%%%%%%%%%%%%%%%%%%
%%%%%%%%%%%%%%%%%%%%%%%%%%%%%%%%%%%%%%%%%%%%%%%%%%%%%%%%%%%%%%%%%%%%%%%%

\medskip

{\bf Acknowledgements.}~Partial support from the KBN
grant 2/P03A/002/24, the POLONIUM project \'EGIDE--KBN
No.~05643SE, and the EU contract HYKE No. HPRN-CT-2002-00282, is
gratefully ackowledged.


\medskip

\begin{thebibliography}{GG}

%\def\dr#1{{\bf [#1]}}
\def\dr#1{}

\vspace{-1mm}


\bibitem{BKL}
S. Benachour, G. Karch \& Ph. Lauren\c cot,
{\it Asymptotic profiles of solutions 
to viscous Hamilton-Jacobi equations}, (submitted).

\vspace{-6pt}

\bibitem{BKLxx} S.~Benachour, H.~Koch \& Ph.~Lauren\c cot, {\it Very
singular solutions to a nonlinear parabolic equation with absorption. II
-- Uniqueness},  Proc. Roy. 
Soc. Edinburgh Sect.~A, to appear.

\vspace{-6pt}

\bibitem{BL99}
S. Benachour \& Ph. Lauren\c cot, 
{\it Global solutions to viscous 
Hamilton--Jacobi equations with irregular initial data}, 
Comm. Partial Differential Equations {\bf 24} (1999), 1999--2021.

\vspace{-6pt}

\bibitem{BL01}
S. Benachour \& Ph. Lauren\c cot, 
{\it Very singular solutions to a nonlinear parabolic equation with 
absorption I. Existence}, Proc. Roy. Soc. Edinburgh Sect. A {\bf
131} (2001), 27--44. 

\vspace{-6pt}

\bibitem{BK99}
M. Ben-Artzi \& H. Koch, 
{\it Decay of mass for a semilinear parabolic equation}, 
Comm. Partial Differential Equations {\bf 24} (1999), 869--881.

\vspace{-6pt}

\bibitem{BSW02}
M. Ben-Artzi, Ph. Souplet \&  F.B. Weissler, 
{\it The local theory for viscous Hamilton-Jacobi equations in Lebesgue 
spaces}, J. Math. Pures Appl. {\bf 81} (2002), 343--378. 

\vspace{-6pt}

\bibitem{EZ} 
M. Escobedo \& E. Zuazua, {\it Large time behavior for
convection-diffusion equations in $\RR^N$,}
 J. Funct. Anal. {\bf 100}
(1991), 119--161. 

\vspace{-6pt}

\bibitem{EVZ1}
M. Escobedo, J.L. V{\'a}zquez \& E. Zuazua, {\it Asymptotic behavior
and source-type solutions for a diffusion-convection equation}, Arch.
Rational Mech. Anal. {\bf 124} (1993), 43--65.

\vspace{-6pt}

\bibitem{FL}
E. Feireisl \& Ph. Lauren\c cot, 
{\it The $L\sp 1$-stability of constant states of degenerate
convection-diffusion equations}, 
Asymptot. Anal. {\bf 19} (1999), 267--288.  

\vspace{-6pt}
 
\bibitem{GGK03}
B.~Gilding, M.~Guedda \& R.~Kersner, \textit{The Cauchy problem for
  $u_t=\Delta u + \vert\nabla u\vert^q$},
J. Math. Anal. Appl. \textbf{284} (2003), 733--755.

\vspace{-6pt}

\bibitem{KS}
G. Karch \& M.E. Schonbek, {\it On zero mass solutions
 of viscous conservation laws,} Comm. Partial Differential Equations {\bf
 27} (2002), 2071--2100.

\vspace{-6pt}

\bibitem{kim}
Y.J.~Kim, {\it An Oleinik type estimate for a convection-diffusion
equation and convergence to $N$-waves}, J. Differential Equations, to
appear. 

\vspace{-6pt}

\bibitem{LS03}
Ph. Lauren\c cot \& Ph. Souplet {\it On the growth of mass for a viscous
Hamilton-Jacobi equation}, J. Anal. Math. {\bf 89} (2003), 367--383.

\vspace{-6pt}

\bibitem{Sm}
J.~Smoller,  {\it Shock waves and reaction-diffusion equations},
Grundlehren der Mathematischen Wissenschaften {\bf 258},
Springer-Verlag, New York-Berlin, 1983.
  
\end{thebibliography}
\end{document}